\documentclass[12pt]{article}
\newcommand{\R}{{\mathbb R}}

\usepackage{graphicx}
\usepackage{amssymb}
\usepackage{amsmath}
\usepackage{amsthm}
\usepackage{color}
\numberwithin{equation}{section}
\def\grad{\nabla}

\newtheorem{theorem}{Theorem}
\newtheorem{lemma}[theorem]{Lemma}

\newtheorem*{theorem*}{Theorem}

\title{ABP and global H\"older estimates for \\
fully nonlinear elliptic equations in unbounded domains
\vskip-0.75cm}
\author{I. Birindelli, I. Capuzzo Dolcetta, A. Vitolo}
\date{}
\begin{document}
\maketitle

\section{Introduction}\label{INTRO}
This paper is focused on two goals about
fully nonlinear degenerate elliptic equations in unbounded domains of $\Omega\subset\R^n$. The first one is to provide some generalization of a well-known result of X. Cabr\'e who proved in \cite{Cab} 
that the Alexandrov-Bakelman-Pucci (ABP, in short) holds  for linear second order uniformly elliptic operators under the following measure theoretic condition on $\Omega$ :
{\it there exist positive real numbers $\sigma$ and $R_0$ 
such that for each $y \in \Omega$ there is a ball $B_{r_y}$ of radius $r_y\le R_0$ such that
\begin{equation*}
y \in B_{r_y}\quad \mbox{and} \quad |B_{r_y}\backslash\Omega_{y}| \ge \sigma|B_{r_y}| \eqno{\rm (G_y)}
\end{equation*}}
where $|\ \, |$ denotes the $n-$dimensional Lebesgue measure and $\Omega_{y}$ is the (connected) component of $\Omega \cap B_{r_y}$ containing $y$.  We will refer to this condition as to (G$_\Omega$) and to $\Omega$ satisfying this conditions as to a cylindrical domain. The
reader will notice a slight simplification with respect to 
the original condition \cite{Cab}, see also \cite{V-Sin} in this respect.

\noindent The second goal is to  prove global H\"older estimates for the same class of operators.
This will require that the unbounded domain satisfies the slightly stronger condition (G$^d_\Omega$), namely :
{\it there exist real numbers $0<\sigma <1$, $d_0>0$ and $K_0>1$ such that, for each $y \in \Omega$, 
condition $(G_y)$ above holds for radii} $$r_y\le K_0d(y)\,\eqno{\rm (d_y)}\,$$ 
and $d(y) \le d_0$. 
%\begin{equation*}
%{\rm dist}(y,\partial\Omega)\equiv d(y)\le d_0. 
%\end{equation*}
\noindent Here, $d(y)$ denotes the distance from point $y$ to $\partial\Omega$.
As mentioned above, this is a slightly stronger condition than  (G$_\Omega$): indeed  (G$^d_\Omega$) implies (G$_\Omega$) 
with $R_0 \le K_0d_0$. 
Observe that the converse is not true. 
In fact, it is immediate to realize that, if we set $R_k^\pm=\{x=(x_1,x_2) \in \mathbb R^2 : \pm x_1 \ge 1,\,x_2=k\}$ and
 $$C=\displaystyle\bigcup_{k \in \mathbb Z}(R_k^+\cup R_k^-),$$ 
then  $\Omega =  \mathbb R^2 \backslash C$ satisfies condition (G$_\Omega$) but not (G$^d_\Omega$).

\noindent In \cite{CS} S. Cho and M. Safonov have proved global H\"older estimates for solutions of the Dirichlet problem, with homogeneous boundary condition,  for linear second-order uniformly elliptic equations in non-divergence form. Their assumption on $\partial \Omega$ is the following `exterior measure' condition:\\
there exists a constant $\sigma>0$ such that for each $y \in \partial \Omega$ and $r>0$ 
\begin{equation*}
|B_{r}(y)\backslash\Omega| \ge \sigma|B_r| \eqno{\rm (A_y)}
\end{equation*}
where $B_r(y)$ is the ball of radius $r >0$, centered at $y$.
\noindent We now describe more precisely the class of operators that will be considered in the present paper. 
Let ${\cal S}_n$ be the space of $n \times n$ real symmetric matrices and recall the definition of the maximal and minimal Pucci operators, acting on ${\cal S}^n$ as 
\begin{align*}
{\cal P}^+_{\lambda,\Lambda}(X)=& \sup_{\lambda I \le A \le \Lambda I}{\rm Tr}(AX) \equiv \Lambda {\rm Tr}(X^+)-\lambda {\rm Tr}(X^-),\\
{\cal P}^-_{\lambda,\Lambda}(X)=& \inf_{\lambda I \le A \le \Lambda I}{\rm Tr}(AX) \equiv \lambda {\rm Tr}(X^+)-\Lambda {\rm Tr}(X^-)\,
\end{align*}
where  $\Lambda \ge \lambda>0$, $I$ is the $n$-dimensional identity matrix and $\rm Tr(X)$ is the trace of $X$.

\noindent The structural condition that we assume on the scalar mapping $F$ defined on  
$\Omega\times\R\times\R^n\times {\cal S}^n$:
{\it there exists $\delta\in(0,1)$ such that \begin{equation}\label{largegrad-subsoln}
 {\cal P}^-_{\lambda,\Lambda}(X)-b^-(x)|\xi|\le F(x,u,\xi,X) \le {\cal P}^+_{\lambda,\Lambda}(X)+b^+(x)|\xi|\;, \ \ 
\end{equation}
for all $x \in \Omega$, $u \ge 0$, $X \in {\cal S}^n$\, and all $\xi$ such that $\delta \le |\xi| \le \frac1\delta$,
and $b^\pm(x)=\max(\pm b(x),0)$ are continuous functions.}

\noindent The above condition has been introduced by C. Imbert in \cite{Imb}, where he proved the ABP maximum principle and Harnack inequality in bounded domains.
If condition (\ref{largegrad-subsoln}) holds with $\delta=0$, then $F$ is in between two uniformly elliptic operators, and therefore we are essentially in the framework of uniform ellipticity.  
On the other hand, should (\ref{largegrad-subsoln}) hold with $\delta>0$, the uniform elliptic bound from above and from below is required only when the gradient term is away from zero and from infinity. As an example , consider
\begin{equation}\label{hom-grad}
F(\xi,X)=|\xi|^\beta\,{\cal P}^+_{\lambda,\Lambda}(X), \, \mbox{with }\, \beta\geq 0,
\end{equation}
or, more generally,
$G(x,\xi,X)=|\xi|^\beta F(x,X)$ with $F$ uniformly elliptic.

% this condition was introduced by C. Imbert in \cite{Imb} . 
\noindent Consider now the Dirichlet problem
\begin{equation}\label{DP}
\left\{\begin{array} {cl}
F(x,u,Du,D^2u)=f(x) & in \ \Omega\\
u=0 & on \ \partial\Omega,
\end{array}\right.
\end{equation}
The first result of this paper is the following uniform estimate of (ABP) type for viscosity solutions of (\ref{DP}).
\begin{theorem}\label{intro1} Assume that $\Omega$ and $F$ satisfy, respectively conditions $(G^d_\Omega)$ and $(\ref{largegrad-subsoln})$  with  $b \in C(\Omega)$ such that $|b(x)| \le b_0$ and $f \in C(\Omega)$. Let $u\in C(\Omega)$
be a bounded viscosity solution of (\ref{DP}).
Then there exists a positive constant $\alpha'=\alpha'(n,\lambda,\Lambda,K_0,b_0d_0,\sigma) \in (0,1)$  such that 
\begin{equation}\label{ABP-dist-|u|}
\sup_{\Omega}|u(y)|d^{-\alpha}(y) \le C_0K_0d_0^\alpha\,\delta_{\alpha,f}
\end{equation}
where 
\begin{equation}\label{f-norm}
\delta_{\alpha,f}=\max\left[d_0^{1-\alpha}\delta\,; \sup_{y \in\Omega} d^{1-\alpha}(y) \|f\|_{L^n(\Omega\cap B_{r_y})}\right]
\end{equation}
and $B_{r_y}$ is a ball provided by condition $(G_y)$. 
%$\delta_{\alpha,f}=\max\,(\delta_{\alpha,f}^+,\delta_{\alpha,f}^-)$ and $C_0 \in \mathbb R_+$ depends 
\end{theorem}
%\noindent This is a straightforward consequence of Theorem \ref{dgamma-thm} and of Corollary \ref{dgamma-soln-thm}.
%COMMENTARE E CITARE RISULTATI PRECEDENTI SU ABP:\\
\noindent For $\alpha=0$ the above result generalizes the ABP estimate by Cabr\'e \cite{Cab} concerning linear uniformly elliptic operators in cylindrical domains. The same result has been proved by Cafagna and Vitolo in \cite{CV} in more general domains satisfying condition (G$_\Omega$) with no condition on the radii $r_y$, for pure second order operators, and in \cite{V-JDEQ} with unbounded radii $r_y$ of at most linear growth, for operators with lower order terms. In the uniformly elliptic fully nonlinear case, the ABP estimate in unbounded domains have been proved first by Capuzzo Dolcetta, Leoni and Vitolo in \cite{CDLV}.

\noindent  As a consequence of Theorem \ref{intro1}, using known interior regularity results (see \cite{CafCab},\cite{KT},\cite{ARV}),  we can derive the global H\"older continuity of viscosity solutions of the Dirichlet problem (\ref{DP}) in a domain satisfying property $(G^d_\Omega)$.
Consider, for $\alpha \in(0,1)$, the H\"older norm
$$
\|u\|_{\alpha,\Omega}=\|u\|_{L^\infty(\Omega)}+[u]_{\alpha,\Omega}\,\,,
$$ where $
[u]_{\alpha,\Omega} = \sup_{x,y \in \Omega}\frac{|u(x)-u(y)|}{|x-y|^\alpha}\,.
$
\begin{theorem}\label{Holder-thm} Assume that $\Omega$ satisfies condition $(G^d_\Omega)$ for some positive real numbers $\sigma <1$, $K_0>1$ and $d_0$, and that $F$ satisfies $(\ref{largegrad-subsoln})$  with  $b \in C(\Omega)$ such that $|b(x)| \le b_0$. Assume, moreover, that $f \in C(\Omega)$ is such that, for some $\overline \alpha >0$,
\begin{equation}\label{delta-bounded}
\sup_{y \in \Omega} d^{1-\overline\alpha}(y)\|f\|_{L^n(\Omega\cap B_{r_y})} < \infty\,.
\end{equation}
Let $u \in C(\overline \Omega)$ be a viscosity solution of  (\ref{DP}) such that $|u(x)|\le M$.
Then there exists $\alpha''=\alpha''(n,\lambda,\Lambda,K_0,b_0,d_0,\sigma)\in(0,1)$ such that
%$u$ is globally H\"older continuous in $\Omega$ and
\begin{equation}\label{Holder-global}
\|u\|_{\alpha,\Omega} \le C\,\delta_{\alpha,f}
\end{equation}
for some $\alpha \in (0,\alpha'']$ where $C$ is a positive constant depending on $n,\lambda,\Lambda,K_0,b_0,d_0,\sigma,\alpha$ and $\delta_{\alpha,f}$ is as in $(\ref{f-norm})$.
\end{theorem}

\noindent Operators like the one in (\ref{hom-grad}) fit into our framework. Nonetheless, the results obtained here in the case $\beta \in (0,1)$ can be generalized to the singular case $-\beta$ with $\beta \in (0,1)$, using the fact that solutions of singular equations are solutions of uniformly elliptic equations with lower order terms. In order to be more precise on this point, consider the equation
\begin{equation}\label{is}
G(x,\grad u,D^2u):=|\grad u|^{-\beta}F(x,D^2u)=f(x)\ \mbox{in}\ \Omega, \mbox{with }\, \beta\in (0,1).
\end{equation}
According to Birindelli-Demengel \cite{BD1} a lower semicontinuous function
$u$ is a supersolution of (\ref{is}) with $\beta\in (0,1)$
if, for each $x_o\in\Omega$, 
either there exists an open ball $B(x_o,\delta) \subset \Omega$, $\delta>0$,
on which  $u$ is constant and 
$f \ge 0$, \\
or for all $\varphi\in {\mathcal C}^2(\Omega)$ such that
$u-\varphi$ has a local minimum at $x_o$ and $\grad\varphi(x_o)\neq
0$ the following inequality holds
$$|\grad \varphi(x_o)|^{-\beta}F(x_o,D^2\varphi(x_o))\leq f(x_o).$$

\noindent Similar  definition can be given for subsolutions. It can be proved, see Lemma \ref{equivalence} in the next Section, that if $u$ is a solution of (\ref{is}) in the sense of the above definition then it is a solution of
$$F(x,D^2u)-|\grad u|^{\beta}f(x)=0\ \mbox{in}\ \Omega$$
in the standard viscosity sense, see \cite{CIL}.

%, while,  at the same time,  it satisfies
% $$\tilde G(x,\xi,M)=F(x,M)-f(x)|\xi|^\beta\geq  {\cal P}^-_{\lambda,\Lambda}(X)-c(x)|\xi|$$
% with .

\noindent We recall that for the whole range $\beta \in (-1,1)$, regularity results and Harnack inequalities have been proved  by Birindelli and Demengel in \cite{BD1,BD2} in the case of bounded smooth domains. See also \cite{BD3} for domains where the boundary may contain conical points.

\noindent As already said, global H\"older results for solutions of second order linear uniformly elliptic equations in bounded domains have been proved by Cho and Safonov \cite{CS}. Interior H\"older continuity estimates in the fully nonlinear setting are due to Trudinger \cite{Tr} and Caffarelli \cite{Caf}; see also \cite{CafCab}, where an extension up to the boundary for regular bounded domains is reported on. Further local results have been proved by \'Swiech \cite{Sw} and by Sirakov \cite{Si}, where also a global $C^\alpha$-regularity result is proved in bounded domains with an uniform exterior cone property. 

\noindent The paper is organized as follows: in the next Section we establish a Growth Lemma for subsolutions, based on a boundary weak Harnack inequality and prove Lemma \ref{equivalence}; in the Section 3 we prove our main result concerning the uniform estimate of the velocity of solutions approaching the boundary value zero and consequently, combining with interior $C^\alpha$-estimates, our global H\"older regularity results.

\section{Preliminaries} \label{preliminar}

We start with a growth lemma (see \cite{CS}). The argument of the proof goes back to \cite{Cab} and can be also found in \cite{CV}, \cite{V-JDEQ} for the linear case, and in \cite{CDLV} considering fully nonlinear operators. The basic tool is the Krylov-Safonov weak Harnack inequality proved by C. Imbert in Theorem 2 of \cite{Imb}, which we will  use in the rescaled version:
\begin{equation}\label{wH-ineq}
\left(\frac1{|B_r|}\int_{B_r}(v)^{p_0}\,dx\right)^{1/p_0}\le C(\inf_{B_r}v +r\max(\delta,\|g^+\|_{L^n(B_{\frac r \tau})}))
\end{equation}
 for positive constants $p_0$ and $C>1$ depending on $n,\lambda,\Lambda,b_0r,\tau$. Here $v$ is a non-negative solution of a second-order degenerate elliptic equation 
\begin{equation}\label{G-supsoln}
G(x,v,Dv,D^2v) \le g(x) 
\end{equation} 
in $B_{\frac r\tau}$, where 
\begin{equation}\label{largegrad-G}
{\cal P}^-_{\lambda,\Lambda}(X)-b^-(x)|\xi| \le G(x,t,\xi,X) \ \  {\rm for \ all} \ \ \xi \in \mathbb R^n :  \ |\xi| \ge \delta,
\end{equation}
for some $\delta \in \mathbb R_+$ and for all $x \in \Omega$, $u \ge 0$, $X \in {\cal S}_n$. We also notice that $b^-(x)$ and $g^+(x)$ are continuous functions.\\
By a standard viscosity argument, if $v$  satisfies (\ref{G-supsoln}) in $D$, then inequality (\ref{wH-ineq}) can be extended up to the boundary (see, for instance, \cite{GT} in the linear case and \cite{CafCab} in the fully nonlinear viscosity setting) introducing the supersolution (defined in all $B_{\frac r\tau}$) 
$$
v^-_m(x)=\left\{\begin{array}{cl}
\inf(v(x),m) & if \ x \in B_{\frac r\tau} \cap D\\
m & if \ x \in B_{\frac r\tau} \backslash \,D,
\end{array}\right.
$$
where $\displaystyle m=\inf_{\partial D \cap B_{\frac r\tau}}v$. In this case (\ref{wH-ineq}) yields
\begin{equation}\label{boundary-wH-ineq}
\left(\frac1{|B_r|}\int_{B_r}(v^-_m)^{p_0}\,dx\right)^{1/p_0}\le C(\inf_{D\cap B_r}v +r\max(\delta,\|g^+\|_{L^n(D\cap B_{\frac r \tau})}))
\end{equation}
with  $p_0$ and $C$ as before.

\begin{lemma}\label{growth-lemma} %2.1 
Let $D$ be a domain of $\mathbb R^n$ and $B_r$ be a ball such that $B_{r}\cap D\neq\emptyset$ and $|B_r\backslash D| \ge \sigma |B_r|$ for $\sigma\in (0,1)$. Let also $B_{\frac r\tau}$ be the concentric ball of radius $\frac r\tau$ with $\tau \in (0,1)$. Let also $\delta \ge 0$ and $F$ be such that the right-hand side of the structure condition $(\ref{largegrad-subsoln})$ holds true. 
Suppose furthermore that $b^+(x)$ and $f^-(x)$ are continuous functions such that
\begin{equation}\label{b-norm}
\sup_{\Omega}\,b^+(x) \le b_0 < \infty.
\end{equation}
%and 
%\begin{equation}\label{f-norm-}
%\delta_{\alpha,f}^-=\max\left(d_0^{1-\alpha}\delta,\, \sup_{y \in\Omega} d^{1-\alpha}(y) \|f^-\|_{L^n(\Omega\cap B_{r_y})}\right).
%\end{equation}
 If $u \in usc(\overline D)$ is a viscosity subsolution, bounded above, of equation
\begin{equation}\label{eqn:degenerate}
F(x,u,Du,D^2u) = f(x)
\end{equation} 
in $D$. Then there exists a positive constant $\theta_0=\theta_0(n,\lambda,\Lambda,b_0r,\sigma,\tau)<1$ such that
\begin{equation}\label{G-lemma-ineq}
\sup_{\Omega\cap B_{r}}u \le \theta_0 \sup_{D\cap B_{\frac{r}\tau}}u^+ + (1-\theta_0)\,\sup_{\partial D \cap B_{\frac{r}\tau}}u^+ + r\,\delta_{f}^-(D\cap B_{\frac{r}\tau})\,,
\end{equation} 
where
\begin{equation}\label{f-norm-local-}
\delta_{f}^-(\Omega)=\max\left(\delta,\|f^-\|_{L^n(\Omega)}\right).
\end{equation}
\end{lemma} 
\noindent{\bf Proof.}  
Passing to $u^+$, which is in turn a viscosity solution of the differential inequality 
$$
F(x,0,Du^+,D^2u^+) \ge -f^-(x),
$$
we set $M_r\equiv \sup_{D \cap B_r}u^+$ and observe that
$$
v(x) = M_{\frac{r}\tau}-u^+(x)
$$
satisfies the differential inequality (\ref{G-supsoln}) with 
$$G(x,t,\xi,X)= -F(x,0,-\xi,-X), \quad g(x)=f^-(x).$$
Then we apply the above inequality (\ref{boundary-wH-ineq}) noting that
\begin{align*}
 m&=\,\inf_{\partial D \cap B_{\frac r\tau}}(M_{\frac{r}\tau}-u^+)\ge M_{\frac{r}\tau} -\, \sup_{\partial D \cap B_{\frac r\tau}}u^+,\\
 \inf_{D\cap B_{r}}v & =\, \inf_{D\cap B_{r}}(M_{\frac{r}\tau}-u^+)=\,M_{\frac{r}\tau}-M_{r}\,.
\end{align*}
We get
\begin{equation*}
\begin{split}
\sigma^{\frac1{p_0}}\,(M_{\frac{r}\tau} - \sup_{\partial D \cap B_{\frac {r}\tau}}u^+) & \le \left(\frac1{|B_{r}|}\int_{B_{r}}(v^-_m)^{p_0}\,dx\right)^{1/p_0}\\
&\le\,  C\left(\inf_{D\cap B_{r}}v +r\,\delta_{f}^-(D\cap B_{\frac{r}\tau})\right)\\
&\le\,  C\left(M_{\frac{r}\tau}-M_{r} +r\,\delta_{f}^-(D\cap B_{\frac{r}\tau})\right)
\end{split}
\end{equation*}
from which
\begin{equation*}
M_{r} \le (1-\tfrac{\sigma^{\frac1{p_0}}}{C})M_{\frac{r}\tau}+ \tfrac{\sigma^{\frac1{p_0}}}{C}\sup_{\partial D \cap B_{\frac {r}\tau}}u^+ + r\,\delta_{f}^-(D\cap B_{\frac{r}\tau})
\end{equation*}
and therefore the assert follows with $\theta_0=1-\tfrac{\sigma^{\frac1{p_0}}}{C}$\,.\hfill$\Box$

\medskip
We end this section with the following lemma concerning solutions of singular equations:
\begin{lemma}\label{equivalence}
If $u$ is a solution of (\ref{is}) with $\beta \in [0,1)$ then it is a solution of
$$F(x,D^2u)-|\grad u|^{\beta}f(x)=0\ \mbox{in}\ \Omega$$
in the standard viscosity sense.
\end{lemma}
\noindent Note that, for $\beta \in [0,1)$, the operator
 $\tilde G(x,\xi,M)=F(x,M)-f(x)|\xi|^\beta$ does satisfy (\ref{largegrad-subsoln}) for $|\xi|\geq\delta$,
 with $b^\pm(x)=f^\mp(x)\delta^{\beta-1}$\,.

\noindent{\bf Proof.}
Let $u$ be a super solution of (\ref{is}) and let  $\varphi\in {\mathcal C}^2(\Omega)$, such that
$u(x)\geq\varphi(x)$ and $u(x_o)=\varphi(x_o)$.
If $\grad\varphi(x_o)\neq 0$, there is clearly nothing to prove.
So we will suppose that $\grad\varphi(x_o)=0$, and for simplicity we will suppose
that $x_o=0$ and $u(0)=\varphi(0)=0$

Without loss of generality we will take $\varphi(x)=\frac{1}{2}\langle Ax,x\rangle$ and suppose that
$u(x)>\varphi(x)$ in a neighbourhood of $0$. We want to prove that
$$F(0,A)\leq f(0)|\grad\varphi(0)|^{\beta}=0.$$
We suppose by contradiction that $F(0,A)>0$. By the ellipticity 
hypothesis on $F$ this implies that $V^+$, the space of eigenvectors
corresponding to positive eigenvalues of $A$, has at least dimension
one. Let $e\in V^+$ be a unitary eigenvector for $A$.

And for $\varepsilon>0$ we introduce $\psi(x)=\varphi(x)+\varepsilon|\langle x,e\rangle |$. 
Let $x_o^\varepsilon\in B_r(0)$ such that 

$$u(x_o^\varepsilon)-\psi(x_o^\varepsilon)=\inf_{x\in B_r(0)}(u(x)-\psi(x)).$$
Observe first that, for  $\varepsilon$ sufficiently small, the minimum is achieved inside i.e. $|x_o^\varepsilon|<r$.
Indeed $u(0)-\psi(0)=0$, while,  for $0<k:=\min_{x\in \partial B_r(0)}(u(x)-\phi(x))$ and  $\varepsilon<\frac{K}{r}$, 
$$\min_{x\in \partial B_r(0)}(u(x)-\psi(x))\geq K-\varepsilon r>0.$$
Remark  also that $\langle x_o^\varepsilon,e\rangle\neq 0$. Indeed suppose that it is zero, then, 
we would have that $\tilde \psi(x)=\varphi(x)+\varepsilon\langle x,e\rangle+u(x_o^\varepsilon)-\psi(x_o^\varepsilon)$
$$u(x)\geq \psi(x)+u(x_o^\varepsilon)-\psi(x_o^\varepsilon)\geq \tilde\psi(x),\ u(x_o^\varepsilon)=\tilde\psi(x_o^\varepsilon).$$ 
Furthermore, $\grad\tilde\psi(x_o^\varepsilon)=Ax_o^\varepsilon+\varepsilon e\neq 0$ since
$$ \langle Ax_o^\varepsilon+\varepsilon e,e\rangle =\varepsilon.$$
So from the equation we get 
$$F(x_o^\varepsilon,A)\leq f(x_o^\varepsilon)|Ax_o^\varepsilon+\varepsilon e|^{\beta}.$$
Observe that for $\varepsilon\rightarrow 0$, $x_o^\varepsilon\rightarrow 0$, so passing to the limit
we get $F(0,A)\leq 0$ a contradiction.

We have obtained that $\psi$ is a test function for $u$ at $x_o^\varepsilon$, with 
$\grad \psi(x_o^\varepsilon)=Ax_o^\varepsilon+\varepsilon e_\varepsilon e\neq 0$ where 
$e_\varepsilon=\frac{\langle x_o^\varepsilon,e\rangle}{|\langle x_o^\varepsilon,e\rangle|}=\pm1$. 

We choose a sequence $\varepsilon_k$ such that $e_{\varepsilon_k}=1$ and 
$\langle Ax_o^{\varepsilon_k}+\varepsilon_k e,e\rangle=\mu \langle x_o^{\varepsilon_k},e\rangle+\varepsilon_k>0$, $\mu$ being the eigenvalue corresponding to $e$.

Finally we can use $\psi$ as as a test function:

$$F(x_o^{\varepsilon_n},A+\varepsilon_n B)\leq |Ax_o^{\varepsilon_n}+\varepsilon_n|^{\beta}f(x_o^{\varepsilon_n}).$$
Passing to the limit we obtain that
$F(0,A)\leq 0$ a contradiction.

\section{Geometric conditions and boundary inequalities} \label{boundary}

We recall that $y \in \Omega$ satisfies condition $(G_y)$ in $\Omega$ with parameter $\sigma \in (0,1)$ if there exists a ball $B_{r_y}$ of radius $r_y$ such that
\begin{equation*}
y \in B_{r_y}, \quad |B_{r_y}\backslash\Omega_{y}| \ge \sigma|B_{r_y}| 
\end{equation*}
where $\Omega_{y}$ is the (connected) component of $\Omega \cap B_{r_y}$ containing $y$. \\
Here we use Lemma \ref{growth-lemma} to obtain a pointwise estimate for viscosity solutions of second-order uniformly elliptic equations in a point $y \in \Omega$ satisfying condition $(G_y)$ in $\Omega$.

\begin{lemma}\label{G-lemma} %3.1
Let $\Omega$ be a domain of $\mathbb R^n$ and suppose that $y$ satisfies condition $(G_y)$ in $\Omega$ with parameter $\sigma \in (0,1)$. Let also $B_{r_y}$ be a ball of radius $r_y \le d_0$ realizing condition $(G_y)$ for a positive constant $d_0$, and $B_{\frac{r_y}\tau}$ be the concentric ball of radius $\frac{r_y}\tau$ with $\tau \in (0,1)$.\\
 Suppose that $u \in usc(\overline \Omega)$ is a viscosity solution, bounded above, of the differential inequality  
\begin{equation}\label{subeqn:degenerate}
F(x,u,Du,D^2u) \ge f(x)
\end{equation} 
in $\Omega$, where $F$ and $f$ satisfy the assumptions of Lemma \ref{growth-lemma}. 
There exists a positive constant $\theta_1=\theta_1(n,\lambda,\Lambda,b_0d_0,\sigma,\tau)<1$ such that 
\begin{equation}\label{point-ineq}
u(y) \le \theta_1 \sup_{\Omega\cap B_{\frac{r_y}\tau}}u^+ + (1-\theta_1)\,\sup_{\partial \Omega \cap B_{\frac{r_y}\tau}}u^+ + 2\,r_y\, \delta_{f}^-(\Omega\cap B_{\frac{r_y}\tau})\,.
\end{equation}
where $\delta_{f}^-(\Omega\cap B_{\frac{r_y}\tau})= \max\left(\delta,\|f^-\|_{L^n(\Omega\cap B_{\frac{r_y}\tau})}\right)$, according to $(\ref{f-norm-local-})$.
\end{lemma}

{\bf Proof.} We set $\tau_0 = (1-\frac{\sigma}4)^{1/n}$, noticing that $\frac34<\tau_0<1$ and considering two different cases according that $(i)$ $y \in B_{\tau_0 r_y}$  or  $(ii)$ $y \not\in B_{\tau_0 r_y}$.

{\it Case $(i)$}. Recall that $\Omega_{y}$ is the component of $B_{r_y}$ containing $y$. Since 
\begin{align*}
|B_{\tau_0 r_y}\backslash \Omega_{y}| \ge  &\,|B_{r_y}\backslash \Omega_{y}|-|B_{r_y}\backslash B_{\tau_0 r_y}|\\
\ge & \,\sigma |B_{r_y}| - (1-\tau_0^n)|B_{r_y}| \ge \frac34\,\sigma |B_{r_y}|
\end{align*}
Therefore we can apply Lemma \ref{growth-lemma} in $D=\Omega_y$ with $r=\tau_0 r_y$ and $\tau=\tau_0$ to get 
\begin{equation}\label{G-lemma-eq1}
\begin{split}
\sup_{x \in \Omega_y\cap B_{\tau_0r_y}}u(x) \le &\, \theta_0 \sup_{\Omega_y\cap B_{r_y}}u^+ + (1-\theta_0)\,\sup_{\partial (\Omega_y \cap B_{r_y})}u^+ + \tau_0r_y\, \delta_{f}^-(\Omega\cap B_{r_y})\\
\le &\, \theta_0 \sup_{\Omega\cap B_{r_y}}u^+ + (1-\theta_0)\,\sup_{\partial \Omega \cap B_{r_y}}u^+ + r_y \,\delta_{f}^-(\Omega\cap B_{r_y})
\end{split}
\end{equation}
and this concludes the proof of (\ref{point-ineq}) with $\theta_1=\theta_0$ in Case (A), since $y \in B_{\tau_0 r_y}$.

\begin{figure}\begin{center}
\includegraphics[height=50mm]{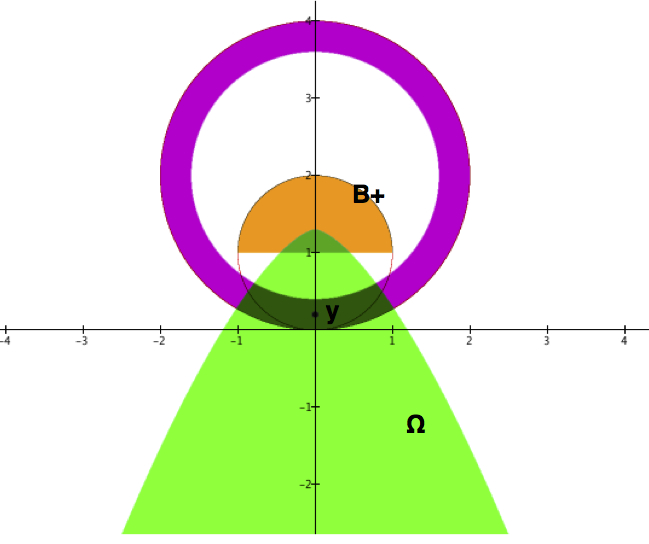}
\end{center}
\caption{Case $(ii)$}
\end{figure}
{\it Case $(ii)$}. Without loss of generality, suppose that $B_{r_y}=B_{r}(P_{r})$ i.e it is the ball of radius $r$ centered at the point $P_{r}=(0,r)$.  In this case  $y \in \Omega\cap \left(B_{r}(P_{r})\backslash B_{\tau_0 r}(P_{r})\right)$ and we may suppose that $y=(0,y_n)$ with $0 \le y_n < (1-\tau_0)r_y$, because $y_n= (1-\tau_0)r$ would mean $y \in \overline B_{\tau_0 r}(P_{r})$ and in this case the result follows from (\ref{G-lemma-eq1}) with $\theta_1=\theta_0$ by continuity.\\ 
Let $y^{o}= (0,\frac{r}2)$. Since $\tau_0>\frac34$, then the upper half-ball
$$B^+=\{x \in B_{\frac{r}2}(y^o)\, /\, x_n \ge \frac{r}2\} \subset B_{\tau_0 r}(P_{r})$$ 
lies in the complement of $\Omega_y \backslash \overline B_{\tau_0 r}(P_{r})$. Consequently, for all $\tau \in (0,1)$,  denoting by $D$ 
the component  of the open set $(\Omega_y \backslash \overline B_{\tau_0 r}(P_r) )\cap B_{\frac r\tau}(y^o)$ containing $y$ , satisfies the measure condition 
$$
|B_{\frac{r}2}(y^o)\backslash D| \ge |B^+| \ge \frac12\,|B_{\frac{r}2}(y^o)|.
$$  
So that we can apply Lemma \ref{growth-lemma} to get
\begin{equation}\label{G-lemma-eq2}
\begin{split}
u(y) \le &\, \theta_\tau \sup_{D\cap B_{\frac{r}{2}}(y^o)}u^+ + (1-\theta_\tau)\,\sup_{\partial D \cap B_{\frac{r}{2}}(y^o)}u^+ + r \delta_{f}^-(D\cap B_{\frac{r}{2}}(y^o))\\
\le &\, \theta_\tau \sup_{\Omega\cap B_{\frac{r_y}\tau}}u^+ + (1-\theta_\tau)\,\max\left(\sup_{\partial \Omega \cap B_{\frac r\tau}(y^o)}u^+,\,\sup_{\Omega_y \cap \partial B_{\tau_0 r_y}}u^+\right) \\
&+ \tfrac12\,r_y \,\delta_f^-(\Omega\cap B_{\frac{r_y}\tau}).
\end{split}
\end{equation} 
If $\displaystyle\sup_{\partial \Omega \cap B_{\frac r\tau}(y^o)}u^+\ge \sup_{\Omega_y \cap \partial B_{\tau_0 r_y}}u^+$, 
we are done and we get (\ref{point-ineq}) with $\theta_1=\theta_\tau$. \\
Otherwise, we again use (\ref{G-lemma-eq1})  obtaining from (\ref{G-lemma-eq2})
\begin{equation}\label{G-lemma-eq3}
\begin{split}
u(y) \le &\, \theta_\tau \sup_{\Omega\cap B_{\frac {r_y}\tau}}u^+ + (1-\theta_\tau)\,(\theta_0 \sup_{\Omega\cap B_{r_y}}u^+ + (1-\theta_0)\,\sup_{\partial \Omega \cap B_{r_y}}u^+)\\
&+\, (1-\theta_\tau)\,r_y \,\delta_f^-(\Omega\cap B_{\frac{r_y}\tau})+\tfrac12\,r_y\, \delta_f^-(\Omega\cap B_{\frac{r_y}\tau})\\
\le &\, (\theta_\tau +(1-\theta_\tau)\,\theta_0)\,\sup_{\Omega\cap B_{\frac {r_y}\tau}}u^+ + (1-\theta_\tau)\,(1-\theta_0)\,\sup_{\partial \Omega \cap B_{r_y}}u^+\\
&+\, \tfrac32\,r_y\, \delta_f^-(\Omega\cap B_{\frac{r_y}\tau})\\
\end{split}
\end{equation} 
which yields (\ref{point-ineq}) with $\theta_1= \theta_\tau +(1-\theta_\tau)\,\theta_0$ and concludes the proof.\hfill$\Box$

\noindent {\it Comment.} Taking the supremum over $y \in \Omega$ in Lemma \ref{G-lemma}, if $\Omega$ satisfies condition $G_\Omega$, being $r_y \le R_0$, we get 
\begin{equation}\label{Cabre-gen}
\sup_{\Omega}u \le \sup_{\partial\Omega}u^+ + CR_0 \max\,(\delta,\sup_{y \in\Omega} \|f^-\|_{L^n(\Omega\cap B_{r_y})})
\end{equation}
for solutions, bounded above, of the differential inequality (\ref{subeqn:degenerate}) in domains $\Omega$ satisfying condition $({G}_\Omega)$. This provides a generalization of the above quoted ABP estimate of Cabr\'e, which corresponds to $\delta=0$.   \\
If $\Omega$ satisfies condition ($G^d_\Omega$), being $r_y \le K_0d(y)$ and $d(y)\le d_0$, from Lemma \ref{G-lemma} we also get for all $\alpha \in (0,1)$ the estimate    
\begin{equation}\label{ABP-modified}
\sup_{\Omega}u \le \sup_{\partial\Omega}u^+ + CK_0d_0^\alpha\,\delta_{\alpha,f}^-
\end{equation}
where $\delta_{\alpha,f}^-$ is given by 
\begin{equation}\label{f-norm-}
\delta_{\alpha,f}^-=\max\left(d_0^{1-\alpha}\delta,\, \sup_{y \in\Omega} d^{1-\alpha}(y) \|f^-\|_{L^n(\Omega\cap B_{r_y})}\right),
\end{equation}
and  $C=C(n,\lambda,\Lambda,b_0,d_0,\sigma,\tau,\alpha)$ is a positive constant.\hfill$\Box$

%\noindent Theorem \ref{dgamma-thm} yields a similar estimate  for $u(x)d^{-\alpha}(x)$.
\begin{theorem}\label{dgamma-thm} 
Let $\Omega$ be a domain of $\mathbb R^n$ satisfying condition $(G^d_\Omega)$ for some positive real numbers $\sigma <1$, $K_0>\max(1,d_0)$. As in Lemma $\ref{growth-lemma}$ we suppose that the right-hand side of the structure condition $(\ref{largegrad-subsoln})$ is satisfied with some $\delta >0$, $b^+(x)$ and $f^-(x)$ are continuous functions, $b^+(x) \le b_0$ in $\Omega$ for a positive real number and recall the definition $(\ref{f-norm-})$,  
\begin{equation*}
\delta_{\alpha,f}^-=\max\left(d_0^{1-\alpha}\delta,\, \sup_{y \in\Omega} d^{1-\alpha}(y) \|f^-\|_{L^n(\Omega\cap B_{r_y})}\right)
\end{equation*}
where $B_{r_y}$ is a ball provided by condition $(G_y)$. \\
Let $u \in usc(\overline \Omega)$ be a viscosity solution, bounded above, of the degenerate elliptic differential inequality (\ref{subeqn:degenerate}).\\
There exists a positive constant $\alpha'=\alpha'(n,\lambda,\Lambda,K_0,b_0d_0,\sigma) \in (0,1)$  such that, if $u\le0$ on $\partial \Omega$, then for $\alpha \in (0,\alpha']$ we have 
\begin{equation}\label{ABP-dist}
\sup_{\Omega}u(y)d^{-\alpha}(y) \le C_0K_0d_0^{\alpha}\,\delta_{\alpha,f}^-
\end{equation}
where $C_0$ is a positive constant depending on $n,\lambda,\Lambda,K_0b_0d_0,\sigma,\alpha$.
\end{theorem}

\noindent {\bf Proof.} Let us consider $\alpha>0$ to be chosen in the sequel (\ref{gamma}) and an arbitrary point $y \in \Omega$, which by assumption satisfies condition $(G_y)$ with parameter $\sigma\in(0,1)$ and a ball $B_{r_y}$ containing $y$  such that $r_y\le K_0\,d(y)$. 

\noindent Since $u$ is a subsolution of equation (\ref{eqn:degenerate}), then
\begin{equation}\label{Pucci-u+}
F(x,u^+,Du^+,D^2u^+)\ge -f^-(x).
\end{equation}
Moreover, since $u$ is supposed to be bounded above, then $u^+(x)(d(x)+\frac1j)^{-\alpha}$ is bounded for all $j \in \mathbb N$, and we set 
$$N_j\equiv \sup_{x \in \Omega}\frac {u^+(x)}{(d(x)+\frac1j)^{\alpha}}<\infty.$$
From inequality (\ref{point-ineq}) of Lemma \ref{G-lemma}, since $u=0$ on $\partial \Omega$ there exists $y^* \in \Omega \cap \overline B_{\frac{r_y}\tau}$ such that
\begin{equation*}\label{}
\begin{split}
u^+(y)\le &\,\theta_1 \sup_{\Omega\cap B_{\frac{r_y}\tau}}u^+ + 2r_y\,\delta_f^-(\Omega\cap B_{\frac{r_y}\tau})\\
= & \,\frac{\theta_1 u^+(y^*)}{(d(y^*)+\frac1j)^{\alpha}}\,(d(y^*)+\tfrac1j)^{\alpha} + 2r_y\,\delta_f^-(\Omega\cap B_{\frac{r_y}\tau})\\
\le & \, \theta_1\,N_j\,(d(y^*)+\tfrac1j)^{\alpha} + 2r_y\, \delta_f^-(\Omega\cap B_{\frac{r_y}\tau})\\
\le &\,  \theta_1\,N_j\,((1+\tfrac1\tau)r_y+\tfrac1j)^{\alpha}+ 2r_y\, \delta_f^-(\Omega\cap B_{\frac{r_y}\tau}),
\end{split}
\end{equation*} 
where $\tau \in (0,1)$ and we have used the fact that $d(x)=dist(x;\partial\Omega) \le (1+\tfrac1\tau)r_y$ for $x \in B_{\frac{r_y}\tau}$ in the last inequality. \\
Since $d(y) \ge \frac {r_y} {K_0}$, we may divide by $(d(y)+\frac1j)^{\alpha}$ the above inequality and get 
\begin{equation*}
\begin{split}
\frac{u^+(y)}{(d(y)+\frac1j)^{\alpha}}\le &\,  \theta_1\, N_j\,\frac{((1+\tfrac1\tau)r_y+\tfrac1j)^{\alpha}}{(\frac{r_y}{K_0}+\frac1j)^{\alpha}}+ \frac{K_0d(y)}{(d(y)+\frac1j)^{\alpha}}\, \delta_f^-(\Omega\cap B_{\frac{r_y}\tau})\\ \\
\le &\,\theta_1\, N_j\,((1+\tfrac1\tau)K_0)^\alpha+ K_0\,d^{1-\alpha}(y)\, \delta_f^-(\Omega\cap B_{\frac{r_y}\tau})\,.
\end{split}
\end{equation*}
Hence, taking the sup over $y \in \Omega$ we obtain
\begin{equation*}
N_j\le N_j\,\theta_1\,((1+\tfrac1\tau)K_0)^\alpha+ K_0d_0^{\alpha}\, \delta_{\alpha,f}^-\,
\end{equation*}
from which, for any positive number
\begin{equation}\label{gamma}
\alpha < \frac{\log\theta_1^{-1}}{\log((1+\tfrac1\tau)K_0)},
\end{equation}
we deduce 
\begin{equation*}
N_j\le C_0K_0d_0^{\alpha}\,\delta_{\alpha,f}^-\,
\end{equation*}
where $C_0=C_0(n,\lambda,\Lambda,b_0K_0d_0,\sigma,\tau,\alpha)$. Finally, letting $j \to \infty$, we get the result.\hfill$\Box$

Suppose that $\Omega$ satisfies condition $(G^d_\Omega)$ as in Theorem \ref{dgamma-thm}. Let 
\begin{equation}
F(x,u,\xi,X) \ge {\cal P}^-_{\lambda,\Lambda}(X)-b^-(x)|\xi|  \ \  for \ all \ \ \xi \in \mathbb R^n :  \ |\xi| \ge \delta,
\end{equation}
and for all $x \in \Omega$, $u \ge 0$, $X \in {\cal S}_n$, where $b^-(x) \le b_0$ is a continuous function.\\
Changing $u$ with $-u$, a similar estimate 
\begin{equation}\label{ABP-dist-inf}
\inf_{\Omega}u(y)d^{-\alpha}(y) \ge -C_0K_0d_0^{\alpha}\,\delta_{\alpha,f}^+,
\end{equation}
where
\begin{equation}\label{f-norm+}
\delta_{\alpha,f}^+ = \max\left(d_0^{1-\alpha}\delta,\,\sup_{y \in \Omega} d^{1-\alpha}(y)\|f^+\|_{L^n(\Omega\cap B_{r_y})}\right),
\end{equation}
can be obtained for a viscosity solution $u \in lsc(\overline \Omega)$, bounded from below and non-negative on $\partial \Omega$, of the uniformly elliptic differential inequality 
\begin{equation}\label{supeqn:degenerate}
F(x,u,Du,D^2u) \le f(x).
\end{equation} 
for a continuous function $f^+(x)$.\\

\noindent {\bf Proof of Theorem \ref{intro1}.} Gathering (\ref{ABP-dist}) and (\ref{ABP-dist-inf}), we get at once Theorem \ref{intro1}.

\noindent {\bf Proof of Theorem \ref{Holder-thm}.} Let $u$ be a continuous viscosity solution of the degenerate elliptic equation  $(\ref{eqn:degenerate})$ in $\Omega$. For the $L^\infty$-norm of $\|u\|_{L^\infty(\Omega)}$ we have already obtained a bound (\ref{ABP-modified}) for all $\alpha\in(0,1)$. \\
To estimate the H\"older seminorm $[u]_{\alpha,\Omega}$, let us take a point $x$ of $\Omega$ and consider the ball $B_r$ of radius $r>0$ centred at $x$. Let
$$
M_r =\sup_{\Omega \cap B_r}u, \quad m_r=\inf_{\Omega \cap B_r}u, \quad \omega(r)=M_r-m_r.
$$ 
We also denote by $\alpha_1$ any exponent $\alpha$ allowed by Theorem \ref{intro1}.

\noindent If $r \ge 1$, note that, since $d(y) \le d_0$ for all $y \in \Omega$, from (\ref{ABP-dist-|u|}) we get
\begin{equation}\label{omega1}
\begin{split}
\omega(r) \le& \,2\,\sup_{y\in\Omega \cap B_r}|u(y)| \\
\le&\,  2d_0^{\alpha_1}\,\sup_{y\in\Omega \cap B_r}d^{-\alpha_1}(y) \,|u(y)| 
\le\, A \delta_{\alpha_1,f}^-r^{\alpha_1},
\end{split}
\end{equation}
where $A=2C_0K_0d_0^{\alpha_1}$.\\

\noindent If $r <1$, we consider separately the cases: a) $r \le \,\frac12\,d(x)$; b) $r > \,\frac12\,d(x)$.\\

\noindent Case a)
Since $r \le \frac12\,d(x)$, then $B_{r/2} \subset B_{r} \Subset \Omega$. In this case we have
$$
M_r =\sup_{B_r}u, \quad m_r=\inf_{B_r}u, \quad \omega(r)=M_r-m_r,
$$
noticing that the non negative functions $u-m_{r}$ and $M_{r}-u$ are solutions of the degenerate elliptic differential inequalities (\ref{subeqn:degenerate}), with $u-m_r$ instead of $u$, and  (\ref{G-supsoln}) with $v=M_r-u$, respectively. \\
By structure assumptions (\ref{largegrad-subsoln}) and (\ref{largegrad-G}), we can apply to $M_r-u$ and $u-m_r$ the Harnack 
inequality of C. Imbert, Corollary 3 of \cite{Imb}, which generalizes the uniformly elliptic case (see for instance \cite{CafCab},
\cite{KT},\cite{ARV}), to obtain
\begin{equation}
\begin{split}
M_{r}-m_{r/2} \le &C(M_{r}-M_{r/2}+r\delta_f^-(B_r) )\\
M_{r/2} - m_{r} \le &C(m_{r/2}-m_{r}+r\delta_f^+(B_r)),
\end{split}
\end{equation}
where 
\begin{equation}\label{f-norm-local}
\delta_f^\pm(B_r) = \max(\delta, \|f^\pm\|_{L^n(B_r)}), \quad \delta_f(B_r) = \max(\delta, \|f\|_{L^n(B_r)})
\end{equation}
and $C$ is a non-negative constant depending on $n,\lambda,\Lambda,b_0K_0d_0$.\\
Adding the two inequalities above, we have
\begin{equation}
\omega(r/2) + \omega(r) \le C(\omega(r)-\omega(r/2)+2r\,\delta_f (B_r) )
\end{equation}
from which we have
\begin{equation}\label{case-a}
\begin{split}
\omega(r/2)&\le \frac{C-1}{C+1}\,\omega(r)+ \frac{2C}{C+1}\,r\,\delta_f (B_r)\\
&\le \gamma\,\omega(r)+ C_1\delta_{\alpha_1,f}\,r^{\alpha_1} \,.
\end{split}
\end{equation}
with $\gamma = \frac{C-1}{C+1}\in (0,1)$, $C_1 =\frac{2C}{C+1} \in \mathbb R_+$ and $\delta_{\alpha_1,f}$ is defined in (\ref{f-norm}) with the exponent $\alpha=\alpha_1$.

\noindent Case b) Here we observe that, if $r>\frac12\,d(x)$, then for $y \in B_{r}$ we have
$$
d(y) \le d(x)+d(x,y)\le 3r,
$$
from which, using (\ref{ABP-dist}) with $\alpha = \alpha_1$, we deduce
\begin{equation}\label{case-b}
\begin{split}
\omega(r) \le &\,2\,\sup_{y\in\Omega \cap B_r}|u(y)| \le 2\,\sup_{y\in\Omega \cap B_r}\left(\frac {3r}{d(y)}\right)^{\alpha_1} \,|u(y)| \\
\le &\,2\cdot 3^{\alpha_1} \left(\sup_{y\in\Omega \cap B_r}d^{-\alpha_1}(y) \,|u(y)|\right)r^{\alpha_1}
\le \,2\cdot 3^{\alpha_1} C_0K_0\delta_{\alpha_1,f} r^{\alpha_1}\,.
\end{split}
\end{equation}

\noindent Letting $B=\max(C_1,2\cdot3^{\alpha_1} C_0K_0)$, from (\ref{case-a}) and (\ref{case-b}) we get 
\begin{equation}\label{case-ab}
\omega(r/2) \le \gamma\,\omega(r)+B \delta_{\alpha_1,f}r^{\alpha_1}\,.
\end{equation} 
for all $r \in \mathbb R_+$.\\
Since $r < 1$, by virtue of (\ref{case-ab}) Lemma 8.23 of \cite{GT}, together with (\ref{omega1}), yields 
\begin{equation}\label{omegar<1}
\begin{split}
\omega(r) \le C_2\left(\omega(1)r^\beta +B \delta_{\alpha_1,f}r^{\mu\alpha_1} \right)\le C_2\delta_{\alpha_1,f}\left(Ar^{\beta} +B r^{\mu\alpha_1} \right)
\end{split}
\end{equation}
for any $\mu \in (0,1)$, where $C_2=C_2(\gamma)$ and $\beta=\beta(\gamma,\mu)$ are positive constants. \\
From (\ref{omega1}) and (\ref{omegar<1}) we therefore obtain
\begin{equation}
\begin{split}
\omega(r) \le C_3 \delta_{\alpha_1,f} r^\alpha
\end{split}
\end{equation}
for $\alpha\le \min(\beta,\mu\alpha_1)$, where $C_3$ is a positive constant depending on $\gamma$, $A$ and $B$.\\
Hence, if $y$ is any other point of $\Omega$, letting $|y-x|=r$, we get 
$$
|u(y)-u(x)|\le  C_3 \delta_{\alpha,f}|x-y|^\alpha
$$
and we are done. \hfill$\Box$


\begin{thebibliography}{30} 
\bibitem {ARV} M.E. Amendola, L. Rossi and A. Vitolo,  Harnack Inequalities and ABP Estimates for Nonlinear Second Order Elliptic Equations in Unbounded Domains, {\it Abstr. Appl. Anal.}(2008).
\bibitem{BD1} I. Birindelli, F. Demengel, Eigenvalue, maximum principle and regularity for fully non linear
homogeneous operators, {\it Comm. Pure and Appl. Analysis},  \textbf{6} (2007), 335-366.
\bibitem{BD2}  I. Birindelli, F. Demengel,  Regularity and uniqueness of the first eigenfunction for singular fully non linear operators, {\em J.  Differential Equations},  \textbf{249}, (2010), 1089-1110.   
\bibitem{BD3}  Birindelli, I.; Demengel, F. Eigenvalue and Dirichlet problem for fully-nonlinear operators in non-smooth domains. {\em J. Math. Anal. Appl. } \textbf{352} (2009), no. 2, 822-835.
\bibitem{Cab} X. Cabr\'e, On the Alexandroff--Bakelman--Pucci estimate and the reverse H\"older inequality for solutions of elliptic and parabolic equations, {\it Comm. Pure Appl. Math.} {\bf 48} (1995), 539�-570.
\bibitem{CV} V. Cafagna and A. Vitolo,  On the maximum principle for second-order elliptic operators in unbounded domains, {\it C. R., Math., Acad. Sci. Paris} {\bf 334} (2002), n.5, 359--363.
\bibitem{Caf} L.A. Caffarelli, Interior a priori estimates for solutions of fully nonlinear equations, {\it Ann. Math.} {\bf 130} (1989), 189--213.
\bibitem{CafCab} L. A. Caffarelli and  X. Cabr\'{e}, ``Fully Nonlinear Elliptic Equations'', AMS Colloquium Publications, Providence, Rhode Island, 1995.
\bibitem{CS} S. Cho and M. Safonov, H\"older regularity of solutions to second-order elliptic equations in nonsmooth domains {\it Boundary Value Problems} \textbf{2007} (2007), Article ID 57928, pp. 24.
\bibitem{CDV} I. Capuzzo Dolcetta and A. Vitolo, A qualitative Phragm\`en�-Lindel\"of theorem for fully nonlinear elliptic equations, {\it J. Differential Equations} {\bf 243} (2007), 578�-592.
\bibitem{CDLV} I. Capuzzo Dolcetta, F. Leoni and A. Vitolo, {The Alexandrov--Bakelman--Pucci weak maximum principle for fully nonlinear equations in unbounded domains}, {\it Commun. Partial Differ. Equations} {\bf 30} (2005), pp. 1863--1881.
\bibitem{CIL} M.G. Crandall, H. Ishii and P.L. Lions, {User's guide to viscosity solutions of second order partial differential equations}, {\it Bull. Am. Math. Soc., New Ser.} {\bf 27} (1992), 1--67.
\bibitem{GT} D. Gilbarg and  N.S. Trudinger, ``Elliptic Partial Differential Equations of Second Order'', $2^{\hbox{\footnotesize nd}}$ ed., Grundlehren Math. Wiss. 224, Springer-Verlag, Berlin-New York, 1983.
\bibitem{Imb} C. Imbert, Alexandroff-Bakelman-Pucci estimate and Harnack inequality for degenerate/singular fully non-linear elliptic equations, {\it Journal of Differential Equations}  {\bf 250} (2011), 1553--1574.
\bibitem{Koi} S.Koike, {\it A Beginners Guide to the Theory of Viscosity Solutions}, MSJ Memoirs 13, Math. Soc. Japan, Tokyo (2004).
\bibitem{KT} S.Koike and Takahashi, Remarks on regularity of viscosity solutions for fully nonlinear uniformly elliptic PDEs with measurable ingredients, {\it Adv. Differential Equations} \textbf{7} (2002), 493--512.
\bibitem{Si} B. Sirakov, Solvability of uniformly elliptic fully nonlinear PDE, {\it Arch. Ration. Mech. Anal.} \textbf{ 195} n.2 (2010), 579--607.
\bibitem{Sw} A. \'Swiech, $W^{1,p}$-interior estimates for solutions of fully nonlinear uniformly elliptic equations, {\it Adv. Differential Equations} {\bf 2} (1997), 1005--1027.
\bibitem{Tr} N.S. Trudinger, Comparison principles and pointwise estimates for viscosity solutions, {\it Rev. Mat. Iberoamericana} {\bf 4} (1988), 453--468.
\bibitem{V-JDEQ} A.Vitolo, On the Maximum Principle for Complete Second-Order Elliptic Operators in General Domains, {\it J. Differ. Equations} {\bf 194} (2003), n.1, pp.166--184. 
\bibitem{Vit} A. Vitolo, On the Phragm\`en�-Lindel\"of principle for second-order elliptic equations, {\it J. Math. Anal. Appl.} {\bf 300} (2004), 244�-259.
\bibitem{V-Sin} A.Vitolo, A Note on the Maximum Principle for Second-Order Elliptic Equations in General Domains, {\it Acta Math. Sin., Engl. Ser.} {\bf 23} (2007), n.11, pp.1955--1966. 
\end{thebibliography}
\end{document}